\pdfoutput=1
\documentclass[12pt]{article}
\usepackage{amsmath,amsfonts,amssymb,color,a4,epsfig,graphics}
\usepackage[latin1]{inputenc}

\newcommand{\N}{{\mathbb{N}}}

\def\preuve{\begin{proof}}

\def\gd{\delta}

\def\go{\omega}

\newtheorem{defi}{Definition}[section]

\newtheorem{lemm}{Lemma}[section]
\newtheorem{prop}{Proposition}[section]
\newtheorem{rem}{Remark}[section]

\newtheorem{theo}{Theorem}[section]

\newtheorem{crit}{Criterion}[section]

\newenvironment{demo}{\noindent {\it Proof.  }
      \begin{quotation}\noindent}{\end{quotation}\hfill$\square $}
\newenvironment{demorefpos}{\noindent {\textbf{Proof of lemma 4.2}}
      \begin{quotation}\noindent}{\end{quotation}\hfill$\square $}
\newenvironment{demorefharnak}{\noindent {\textbf{Proof of lemma 4.1}}
      \begin{quotation}\noindent}{\end{quotation}\hfill$\square $}
\newenvironment{demorefharmo}{\noindent {\textbf{Proof of Theorem 4.1}}
      \begin{quotation}\noindent}{\end{quotation}\hfill$\square $}

\renewcommand{\footnote}[1]{\footnotetext{#1}}

\begin{document}

\footnote{ \textbf{Keywords:} Infinite graph; weighted graph Laplacian; Schr\"{o}dinger operator; essentially self-adjoint.}
\footnote{ \textbf{Math Subject Classification (2000):} 05C63, 05C50, 05C12, 35J10, 47B25.}
\title{Essential self-adjointness \\for combinatorial Schr\"odinger
operators \\I- Metrically complete graphs}
\author{ Nabila Torki-Hamza
\footnote{Universit\'e de Carthage; Facult\'e des Sciences de Bizerte;
Math\'ematiques et Applications (05/UR/15-02);
 7021-Bizerte (Tunisie);
{\tt nabila.torki-hamza@fsb.rnu.tn};
{\tt torki@fourier.ujf-grenoble.fr}}}

\date{}

\maketitle

%%%%%%%%%%%%%%%%%%%%%%%%%%%%%%%%%%%%%%%%%%%%%%%%%%%%%%%%%%%%%%%%%%%%%%%%%%%%%%%%%%%%%%%%%%%%%%%%%%%%%%%%%%%%%%%%%%%%%%%%%%%%%%%%%%%%%%%%%%%%%

\begin{abstract}
We introduce the weighted graph Laplacian $\Delta_{\omega,c}$
 and the notion of Schr\"{o}dinger operator of the form $\Delta_{1,a}+W$ on a
 locally finite graph $G$.\\
 Concerning essential self-adjointness, we extend Wojciechowski's and
 Dodziuk's results for graphs with vertex constant weight.
The main result in this work states that on any metrically complete weighted
 graph with bounded degree, the Laplacian $\Delta_{\omega,c}$  is essentially
 self-adjoint and the same holds for
 Schr\"{o}dinger operators provided the associated quadratic form is bounded from below.
 We construct for the proof a strictly positive and harmonic function which allows us to write
 any Schr\"{o}dinger operator $\Delta_{1,a}+W$ as a Laplacian  $\Delta_{\omega,c}$ modulo
 a unitary transform.

\end{abstract}

\section{Introduction}

This work is the first of a series of three articles (the others are
[5] et [6]) dealing with spectral theory of Laplacians and Schr\"{o}dinger
operators on infinite graphs. \\
We extend for infinite graphs some classical results of Laplacians and Schr\"{o}dinger
operators on non compact Riemannian manifolds.\\
One of the main results, Theorem 6.2, states that the Laplacian
of a metrically complete weighted graph with bounded degree
is essentially self-adjoint.
Theorem 1.3.1 in [25] and Theorem 3.1 in [12] are Corollaries of this result.
The notion of completeness used here for the weighted graphs is related
to a distance built according to the vertex weight and the edge conductance.
The second article discusses the non complete case for which we will give conditions
of potential increasing to insure essential self-adjointness of a
Schr\"{o}dinger operator. And the third article
deals with the case of combinatorial Schr\"{o}dinger operators with magnetic fields.\\
One of the classical famous questions in mathematical physics is to find conditions
of essential self-adjointness for Schr\"{o}dinger operators.
Many works study this problem in the case of ${\mathbb{R}}^{n}$.
It is mentioned in [1], that the first article in this topic is that of Weyl [23], and the classical results
can be found in the famous four-books [16] of Reed-Simon. \\
Later Gaffney proved, in [9] and [10] (see also [3] and [21]), that the Laplacian
of a complete Riemannian manifold is essentially self-adjoint.
And it is proved (see [16], [18] and [19]) that on a complete Riemannian manifold, a Schr\"odinger operator
 is essentially self-adjoint if the potential satisfies a bounded condition.
The Beltrami-Laplacian on a Riemannian manifold has many analogous on graphs: Laplacians on
quantum graphs (see [8], [14], [2]); combinatorial Laplacians (see [4],
[25], [11],[12]), or physical Laplacians (see [24]).\\

In this article, we introduce a different Laplacian, denoted by $\Delta_{\omega,c}$~
"the \emph{weighted graph Laplacian}", for a  locally finite weighted graph with a
 weight $\omega$ on the vertices and a conductance $c$ on the edges.
  It generalizes the ``combinatorial Laplacian'' of [25] (it is nothing but $\Delta_{1,1}$);
  as well as the ``graph Laplacian'' of [12] (which is  $\Delta_{1,c}$)~.
This notion had been introduced in the  case of finite graphs, see [4] and [22].

In Section 2, we give immediate properties of the weighted graph Laplacian
$\Delta_{\omega,c}$ and we prove that it is unitary equivalent to a
Schr\"{o}dinger operator: $\Delta_{1,a}+W$, modulo a diagonal transform.

In Section 3, we prove that  if the weight $\omega$ is constant,
the operator $\Delta_{\omega,c}$ is essentially self-adjoint.
The idea is inspired from Wojciechowski's method [25].
The same procedure can be used to prove that if we add
a potential $W$ which is bounded from below, it stills essentially self-adjoint.

Section 4 consists on building a
strictly positive function $\Phi$ which is harmonic for a positive
Schr\"{o}dinger operator. We need for this construction a discrete version
of local Harnack inequality, a solved Dirichlet problem
and  a discrete minimum principle.

We use such function $\Phi$, in Section 5,
to prove the important result that
\textit{every positive Schr\"{o}dinger operator is unitary
equivalent to a Laplacian}.

 In Section 6, we consider graphs with bounded degree.
For a given Schr\"{o}dinger operator $\Delta_{1,a}+W$~, we introduce a distance
$\delta_{a}$ on the graph, and we prove that
\textit{for a metrically complete weighted graph, the Schr\"{o}dinger
operator $\Delta_{1,a}+W$ is essentially self-adjoint if its
quadratic form is bounded from below.}
 Then we give a counter-example to show that this Theorem is not a particular case of Theorem 3.2.
In the same section, we deduce a similar Theorem for
the Laplacian $\Delta_{\omega,c}$~. It is the main result of this
article which is a generalization for metrically complete graphs
of Gaffney's theorem.\\

This paper is in fact a translation of [21]. We correct the end of the proof of Lemma \ref{Harnack}  and
we add Remark \ref{K-L}. Furthermore we update some recent references.

%%%%%%%%%%%%%%%%%%%%%%%%%%%%%%%%%%%%%%%%%%%%%%%%%%%%%%%%%%%%%%%%%%%%%%%%%%%%%%%%%%%%%%%%%%%%%%%%%%%%%%%%%%%%%%%%%%%%%%%%%%%%%%%%%%%%%%%%%%
\section{Preliminaries}

Let $G$ an infinite locally finite connected graph. We denote by $V$
 the set of vertices and by $E$ the set of edges.
If  $x$ and $y$ are two vertices of $V$~, we denote by $x\sim y$ if they are connected by an edge
which would be indicated by $\{ x,y\} \in E$~.
Sometimes when $G$ is assumed to be oriented,
we denote by $[x,y]$ the edge from the origin $x$ to the extremity $y$~,
and by $\overline{E}$ the set of all oriented edges.
It is to mention that no result depends on the orientation.\\
The simplest natural infinite locally finite connected graph is the \emph{graph $\N$}
 which is the graph with $V\equiv \N$ and $E=\{\{n,n+1\};n\in \N \}$.

The space of real functions on the graph $G$ is considered as the space of real
functions on $V$ and is denoted by
 $$C( V) =\{ f:~V \longrightarrow {\mathbb{R}}\}$$
 and $C_0\left( V\right)$ is its subset of
finite supported functions.
We consider, for any weight $\omega:~V\longrightarrow {\mathbb{R}}_+^{\star}$ on vertices, the space
\begin{equation*}
l_{\omega}^2 \left( V\right)= \{ f:~V\longrightarrow {\mathbb{R}}~;~
\sum_{x\in V}\omega_x^2\vert f \left( x \right) \vert^2\ <
\infty \} ~.
\end{equation*}
It is a Hilbert space when equipped by the scalar product given by
\begin{equation*}
\langle f,g \rangle_{ l_{\omega}^2 } =\sum_{x\in V} \omega_x^2 f\left(
x\right).g\left( x\right)
\end{equation*}
This space is isomorphic to
\begin{equation*}
l^2 \left( V\right)= \left\lbrace f:~V\longrightarrow {\mathbb{R}}~;~\sum_{x\in
V}\vert f \left( x\right) \vert^2\ < \infty\right\rbrace
\end{equation*}
with respect to the unitary transform
\begin{equation*}
U_\omega :~l_{\omega}^2 \left( V\right)\longrightarrow l^2 \left( V\right)
\end{equation*}
defined by
\begin{equation*}
U_\omega\left( f\right) =\omega f~
\end{equation*}
which preserves the set $C_0\left( V\right)$ of finite supported real functions on $V.$
\begin{rem}
If the weight $\omega$ is constant equal to $\omega_0 \ >0$ (ie. for
any vertex $x\in V$, we have $\omega_x=\omega_0$)~, then
\begin{equation*}
l_{\omega_0}^2 \left( V\right)=l^2\left( V\right)~.
\end{equation*}
\end{rem}

\begin{defi}
The \emph{weighted graph Laplacian on $G$ } with the  \emph{vertex weight}
$\omega:~V\longrightarrow {\mathbb{R}}_+^{\star}$ and the
\emph{edge conductance} $c:~E\longrightarrow {\mathbb{R}}_+^{\star}$,
 is the operator on $l_{\omega}^2\left( V\right)$~, which is denoted by
$\Delta_{\omega,c}~$ and is given by:
\begin{equation}
\left( \Delta_{\omega,c}f\right)\left( x\right) = \dfrac{1}{\omega_{x}^{2}}
~\sum_{ y\sim x }~c_{ x,y}\left( f\left(
x\right) -f\left(y \right) \right)
\end{equation}
for any function $f$ in $l_{\omega}^2\left( V\right)$ and any vertex $x$ in $V$~.
\end{defi}

\begin{rem}
These Laplacians satisfy elementary properties, some of them are taken from
[4] and [7]~:
\begin{enumerate}
  \item The operator $\Delta_{\omega,c}$ is symmetric on $l_{\omega}^2\left( V\right)$
with domain $C_0(V)$~ and its associated quadratic form given by
\begin{equation*}
Q_c (f)=\sum _{ {\lbrace x,y\rbrace}\in E } c_{x,y}(f(x)-f(y))^2~
\end{equation*}
is  positive.
  \item If the weight $\omega$ is in $l^2\left( V\right)$~, the constant
functions are $\Delta_{\omega,c}-$harmonic.
  \item The sums in the expression
  of $\Delta_{\omega,c}$ are finite as the graph $G$
is locally finite, so this operator is well defined on $C_0(V)$~.
  \item It is a local operator, in the sense that $\left(\Delta_{\omega,c}f
\right)\left( x\right)$ depends only on the  values of $f$ on the neighbors of
$x$. The Laplacian $\Delta_{\omega,c}$ can be considered
as a  \textrm{differential operator} on the graph $G$~.
  \item It is an \textrm{elliptic operator,} as for any edge $%
\lbrace x,y\rbrace$ of the graph $G$~, the coefficient $c_{\lbrace x,y\rbrace}$ does not vanish.
  \item The weight $c$ does not depend on the orientation,
  and we have: $~c_{ x,y}=~c_{ y ,x}~,$
for each neighbors $x$ and $y$~.
\end{enumerate}
\end{rem}
To work on the same function space $l^2\left(V\right)$, we use the unitary transform
$U_\omega$~. More precisely, the Proposition
2.1, asserts that $\Delta_{\omega,c}$ is unitary equivalent to a Schr\"{o}dinger
operator of the graph $G$. Let us first give the definition of a combinatorial
Schr\"odinger operator.
\begin{defi}
A \emph{Schr\"odinger operator of the graph} $G$ is an operator
of the form $\Delta_{1,a} +W $ acting on functions of $l^2\left( V\right)$~,
where the potential $W$ is a real function on $V$
and the conductance $a$ is a strictly positive function on $E$~.
\end{defi}

\begin{prop}
 If
\begin{equation*}
\widehat{\Delta}= U_\omega ~\Delta_{\omega,c}~U_\omega ^{-1}~
\end{equation*}
then $\widehat{\Delta}$ is a Schr\"odinger operator
of $G$ and we have
\begin{equation*}
\widehat{\Delta}=\Delta_{1,a}+W
\end{equation*}
where  $a$ is a strictly positive function on $E$ given by:
\begin{equation*}
a_{x,y}=\dfrac{c_{x,y}}{\omega_x\omega_y}
\end{equation*}
 and the potential $W:V\longrightarrow {\mathbb{R}}$ is given by~:
\begin{equation*}
W=-\dfrac{1}{\omega}~\Delta_{1,a}\omega~.
\end{equation*}
\end{prop}

\begin{demo}
For any $g$ in $ C_0\left( V\right)$ and $x$ in $V $, we have:
\begin{align}
\left( \widehat{\Delta}g\right) \left(x \right)&= \omega_{x}\left(
\Delta_{\omega,c}\mathit{U_\omega}^{-1}g\right) \left( x\right)  \notag \\
&=\dfrac{1}{\omega_{x}}\sum_{ y\sim x}~c_{
x,y}\left( \dfrac{g\left(x \right) } {\omega_{x}}-\dfrac{g\left(y
\right) }{\omega_{y}}\right)  \notag \\
&=\sum_{y\sim x} \dfrac{c_{ x,y}} {%
\omega_{x}\omega_{y}}\left( g\left( x\right)-g\left( y\right)\right)
+g\left( x\right)\dfrac{1}{\omega_{x}}\sum_{y\sim x}
c_{ x,y}\left(\dfrac{1} {\omega_{x}}-\dfrac{1}{\omega_{y}}
\right)  \notag \\
&=\left( \Delta_{1,a} g\right)\left( x\right)+W\left(x \right)g\left(
x\right)~,  \notag
\end{align}
where $\Delta_{1,a}$ denotes the Laplacian on $G$ weighted by the vertex
constant weight $\omega\equiv1$ and by
the strictly positive function $a$ on $E$ given by:
\begin{equation*}
a_{ x,y}=\dfrac{c_{ x,y}}{\omega_{x}\omega_{y}}
\end{equation*}
and where the potential $W~:~V\longrightarrow {\mathbb{R}}$ is given by:
\begin{equation*}
W(x)=\dfrac{1}{\omega_{x}}\sum_{y\sim x} c_{
x,y}\left(\dfrac{1}{\omega_{x}} -\dfrac{1}{\omega_{y}}\right)= -%
\dfrac{1}{\omega_x}\left( \Delta_{1,a}\omega\right) \left( x\right)~.
\end{equation*}
\end{demo}

\begin{rem}
In Proposition 2.1, the function $W$ might be obtained
strictly negative, while the Laplacian is positive~: we can take for example, the
graph $G$ with $V={\mathbb{N}}^{\star}$ and $n\sim n+1$ for any $n\in {%
\mathbb{N}}^{\star}$, and we suppose that $G$ is weighted by the vertex weights
$\omega_{n}=\dfrac{1}{n}$ and the edge conductance
$c_{ n,n+1}=\left( n+1\right)^{2}$; then we find
$W\left( n\right)=-n\left( 2n+1\right)\ <0 ~.$
\end{rem}
%%%%%%%%%%%%%%%%%%%%%%%%%%%%%%%%%%%%%%%%%%%%%%%%%%%%%%%%%%%%%%%%%%%%%%%%%%%%%%%%%%%%%%%%%%%%%%%%%%%%%

\section{Extension of Wojciechowski and Dodziuk's results}

Our first two theorems are extensions of results due to
 J. Dodziuk [7] and R.K. Wojciechowski [25] concerning
essential self-adjointness. Let us remind the following definitions.

\begin{defi}
An unbounded symmetric linear operator on a Hilbert space is called
 \emph{essentially self-adjoint} if it has a unique self-adjoint
extension.
\end{defi}

To prove essential self-adjointness,
we use the following useful and practical criterion,
  deduced from Theorem X.26 in [17]~.

\begin{crit}
The definite positive symmetric operator $%
\Delta:~C_0\left(V \right)\longrightarrow l^2( V ) $ is essentially
self-adjoint if and only if $Ker\left(
\Delta^{\star}+1\right)=\lbrace 0\rbrace$~.
\end{crit}

From the definition the adjoint operator $\Delta^{\star}$ of $\Delta~:~C_0\left(V
\right)\longrightarrow l^2\left( V\right)$~, we can deduce:
\begin{equation*}
Dom\left( \Delta^{\star}\right)=  \lbrace f\in l^2\left( V\right)~;~\Delta
f\in l^2\left( V\right)\rbrace~.
\end{equation*}

Using an idea in the proof of Theorem 1.3.1 of [25]~, we prove the following result:

\begin{theo}
If the weight $\omega$ is constant on $V$ then for any
conductance $c$ on $E$, the Laplacian $\Delta_{\omega,c}$~,
with domain $C_0\left(V \right)$~, is essentially
self-adjoint.
\end{theo}

\begin{demo}
Let $\omega_0$ a strictly positive real number, and $\omega\equiv\omega_0$ on
$V~.$ We consider a function $g$ on $V$ satisfying:
\begin{equation*}
\Delta_{\omega_0,c} ~g+g=0~.
\end{equation*}
Let us assume that there is a vertex $x_0$ in $V$ such that $g\left(x_0 \right)\ >0~.$
\newline
The equality
\begin{equation*}
\Delta_{\omega_0,c} ~g\left( x_0\right) +g\left( x_0\right) =0
\end{equation*}
implies
\begin{equation*}
\dfrac{1}{\omega_0^2}\sum_{ y\sim x_0}c_{ x_0,y}\left( g\left(
x_0\right) -g \left( y\right) \right)+g\left( x_0\right) =0~.
\end{equation*}
Then there exists at least a vertex $x_{1}$ for which $g\left( x_0\right) \
< g\left( x_{1}\right)~,$ since $\omega_0\ >0~$ and $c_{ x,y} \
> 0~$ for any edge $\lbrace x,y\rbrace$ in $E~.$  We repeat the procedure with $x_{1}~...$
Hence we build a strictly increasing sequence of strictly positive
real numbers $\left( g\left(x_{n} \right) \right)_{n}~.$
We deduce that the function $g$ is not in $l^2\left( V\right)~.$ \\
A similar way is used to have the same conclusion when we take the assumption $g\left(x_0\right) \ < 0~.$
\end{demo}
\begin{rem}
Theorem 1.3.1 of Wojciechowski in [25] deals with the Laplacian $\Delta_{1,1}~,$
so it is a particular case of Theorem 3.1~.
\end{rem}
We can prove similarly the following Theorem.
\begin{theo}
If $W~:~V\longrightarrow {\mathbb{R}}~$ is a bounded from below potential
 and if $\omega_0$ is a constant weight on $V~,$ then for any
conductance $c:~E\longrightarrow {\mathbb{R_{+}}}~,$
 the Schr\"odinger operator $\Delta_{\omega_0,c}+W~,$
with domain $C_0\left(V\right)~,$ is essentially self-adjoint.
\end{theo}

\begin{demo}
Let $\kappa$ a real number bounding from below the potential $W$~. We proceed as
in the proof of Theorem 3.1~, considering a function $g$ on $V$ satisfying:
\begin{equation*}
\Delta_{\omega_0,c} g+Wg+\kappa_{1}g=0~,
\end{equation*}
avec $\kappa +\kappa_{1} \geq 1~.$
\end{demo}

\begin{rem}
 J. Dodziuk states in Theorem 1.2  (see [7]) that the operator $A+W$ is
essentially self-adjoint when $A$ is a bounded positive symmetric operator
on $l^2(V)$ and when the potential $W$ is bounded from below.\\
 The operator $A$ is $\Delta_{1,c}$ in Theorem 3.2~, and we can conclude that
this Theorem is more general than Dodziuk's,
since the Schr\"{o}dinger operator $\Delta_{1,c} +W$ is essentially self-adjoint
when $W$ is bounded from below, even if the operator $A=\Delta_{1,c}$ is unbounded on $l^2(V)$~,
taking for example the locally finite graph $G$
with unbounded degree and affected of a constant conductance $c\equiv1~.$
\end{rem}
%%%%%%%%%%%%%%%%%%%%%%%%%%%%%%%%%%%%%%%%%%%%%%%%%%%%%%%%%%%%%%%%%%%%%%%%%%%%%%%%%%%%%%%%%%%%%%%%%%%%%%%%%%%%%%%%%%%%%%%%%%%%%%%%%%%%%%%%

\section{Harmonic function on vertices}

We are going to build a function $\Phi$ strictly positive and
harmonic on $V$ which is useful on Section 5~.

\begin{theo}
Let $P$ a Schr\"odinger operator on the graph $G$
 such that
\begin{equation*}
\langle Pf,f\rangle_{l^{2}} > 0~,
\end{equation*}
for any $f$ in $C_0\left( V\right)\setminus\lbrace
0\rbrace ~.$\\
Then there exists a $P-$harmonic strictly positive function $\Phi$
  on $V$~.
\end{theo}

The proof of Theorem 4.1 uses Lemma 4.1 which
gives a local Harnack's inequality for graphs.
At first we present the following definitions:

\begin{defi}
A \emph{subgraph} $G\prime$ of $G$ is a graph having the set of its
vertices included in $V$ and the set of its edges a subset of $E$~.
\end{defi}

\begin{defi}
For a subgraph $G\prime$ of $G $ with the set of vertices $K~,$ we mix up $G\prime$ with $K$, and we define:
\begin{itemize}
\item the \emph{interior} of $K$ denoted by $\overset{\circ}{K}$
\begin{equation*}
\overset{\circ}{K}=\{ x\in K;~y\sim x\Rightarrow y\in K\}
\end{equation*}

\item the \emph{boundary} of $K$ denoted by $\partial K$
\begin{equation*}
\partial K= K\setminus \overset{\circ}{K}= \left\lbrace x\in K~;~\exists ~y\in
V\setminus K,~y\sim x \right\rbrace
\end{equation*}

\item $K $ is \emph{connected} if and only if for any vertices
$x,y$ in $K$, there exist vertices $x_1,x_2,...,x_n$~, such that
\begin{equation*}
x_i \in K,~ x_1=x,~ x_n=y,~ \{x_i,x_{i+1}\}\in E(G\prime)
\end{equation*}
for any $1\leq i\leq n-1$~.
\end{itemize}
\end{defi}

\begin{lemm}\emph{(Harnack)}\label{Harnack}
Let $P$ a Schr\"odinger operator on the graph $G~.$ Let $G\prime$ a
sub-graph of $G$, and let us denote its set of vertices by $K$.
We assume that the interior of $K$ is finite connected.
Then there exists a constant $k\ >0$ such that,
 for any function
$\varphi:~V \longrightarrow {\mathbb{R}} $
strictly positive on $K$ and satisfying
\begin{equation*}
 \ (P\varphi \ )\upharpoonright \overset{\circ}{K} \equiv 0~,
 \end{equation*}
we have:
\begin{equation*}
\dfrac{1}{k} \leq \dfrac{\varphi\left( x\right) }{\varphi\left( y\right)}%
\leq k
\end{equation*}
for any $x,y$ in $\overset{\circ}{K}~.$
\end{lemm}

The resolution of the Dirichlet problem given by Lemma 4.2
is useful to prove Theorem 4.1~.

\begin{lemm}\emph{(Dirichlet)}
Let $P$ a Schr\"odinger operator on the graph $G$
such that for any $f\in C_0\left( V\right)\setminus\lbrace0\rbrace$
we have
\begin{equation*}
\langle Pf,f \rangle_{ l^2 }>0~.
\end{equation*}
Then for any subgraph $G\prime$ of $G$ such that the interior of the
set $K$ of its vertices is finite connected and for any function
$u:\partial K\longrightarrow {\mathbb{R}}~$, there exists a unique
function $f$ on $K$ satisfying the following two conditions:
\begin{description}
  \item[(i)] $\ (Pf\ )\upharpoonright \overset{\circ}{K}\equiv0~.$
  \item[(ii)] $f\upharpoonright \partial K\equiv u~.$
\end{description}
Furthermore, if $u$ is positive and not identically null, then $f$ is strictly
positive on $\overset{\circ}{K}~.$
\end{lemm}

To prove the strict positivity in Lemma 4.2~, we will use a discrete version
of the "minimum principle" , given by Lemma 4.3 in [7]~.

\begin{lemm}\emph{(Minimum principle)}
Let $P=\Delta_{1,a}+W$ Schr\"odinger operator
on the graph $G~,$ where $W\ >0~,$ and let $G\prime$ a subgraph of $G$
such that the set $K$ of its vertices has a finite connected interior.
We assume that there exist a function $f$ satisfying
\begin{equation*}
\langle Pf,f\rangle_{l^{2}} \geq0
\end{equation*}
 in the interior of $K$ and an interior vertex $x_{0}$ so that $f\left( x_{0}\right)$
 is minimum and negative. Then $f $ is constant on $K~.$
\end{lemm}

\begin{demorefpos}
We proceed by steps.
\begin{itemize}
\item For the uniqueness, we suppose the existence of two functions
$f$ and $g$ with finite support in $K$ satisfying the two
conditions of the theorem. Then it follows that
$P\left(f-g\right)\upharpoonright \overset{\circ}{K}\equiv0~,$
and $\left(f-g\right)\upharpoonright \partial K\equiv0~.$
This implies
$$\left\langle P\left( f-g\right),f-g\right\rangle_{l^{2}}=0~.$$
By the hypothesis on $P$, we deduce the nullity of
$\left(f-g\right).$

\item The uniqueness give the existence since the function space on $K$
is finite dimensional.

\item We take a positive not identically null function $u$
and we argue by contradiction, to show that $f$ is strictly positive
 in the interior of $K$~. Assume the existence of a vertex in $\overset{\circ}{K}$
 which has a negative image by $f$.
Let $x_{0}$ the vertex realizing the minimum of $f$ on
$\overset{\circ}{K}$  which is finite and connected.
Thus we have
$f\left(x_{0}\right)\leq0 $ and $Pf$ vanishes on $\overset{\circ}{K}$~.
And from Lemma 4.3~, the function $f$ is constant and negative
on $K$~. This contradicts $f\upharpoonright\partial K\equiv u$
since the function $u $ is supposed a non identically null function. Hence $f$ est strictly
positive on $K~.$
\end{itemize}
\end{demorefpos}

The proof of Harnack's inequality is inspired from
proofs of Lemma 1.6 and Corollary 2.3 in [7]~, noticing that
the constant $k$ does not depend on the function $\varphi~.$

\begin{demorefharnak}
Let us consider a finite subgraph $K$ with a connected interior, and a
function $\varphi :V\longrightarrow {\mathbb{R}}$ strictly positive on
the set of vertices of $K$ and $P-$harmonic on  the set of the vertices
of  $\overset{\circ}{K}$~. Let $x$ and
$y$ two vertices of $\overset{\circ}{K}~.$

\begin{enumerate}
\item[(i)] First we suppose that $\lbrace x,y\rbrace$ is an edge.\\
As $\left(P\varphi \right)\left( x\right)=0$~, ie
\begin{equation*}
\sum_{z\sim x}~a_{ x,z}\left[\varphi\left(
x\right)-\varphi\left( z \right) \right]  +W\left(x \right) \varphi\left(
x\right)=0~,
\end{equation*}
then
\begin{equation*}
\left( \sum_{z\sim x}~a_{ x,z}\right)\varphi\left( x\right)
+W\left(x \right) \varphi\left( x\right) =\sum_{z\sim x}~a_{
x,z} \varphi\left( z \right)~.
\end{equation*}
By the positivity of the functions $\varphi$ and $a$~,We obtain the
 following inequality:
\begin{equation*}
\left[W\left(x \right) +\sum_{z\sim x}~a_{ x,z}\right]
\varphi\left( x\right) \geqslant a_{ x,y} \varphi\left(
y\right)~.
\end{equation*}
Let us denote by ~$\alpha=\min\left\lbrace a_{ r,s} ;r,s\in
K , r\sim s\right\rbrace$ and
\begin{equation*}
A=\sum_{r,s\in K,r\sim s} a_{ r,s }~.
\end{equation*}
The finiteness of $K$ induce $\alpha \ > 0$ et $A\ <\infty$~. Hence denoting:
\begin{equation*}
k_{0}=\dfrac{\max \left( 0,\max_{K} W\right)+A }{\alpha}~,
\end{equation*}
we have: $k_{0}\ >0$~, and we find
\begin{equation*}
\dfrac{1}{k_{0}}\leq \dfrac{\varphi\left( x\right) }{\varphi\left( y\right) }%
\leq k_{0}~.
\end{equation*}

\item[(ii)] If the vertices $x$ and $y$ are not neighbors,
 by the connectedness of $\overset{\circ}{K}$~, there exists a path connecting $x$
to $y$ in $\overset{\circ}{K}$~.
Let us denote the consecutive vertices of the path by: $x_{1}=x,x_{2},x_{3}$,...,$x_{d}=y$.
So we have:
\begin{equation*}
\dfrac{1}{k_{0}}\leq\dfrac{\varphi\left( x_{i}\right) } {\varphi\left(
x_{i+1}\right)}\leq k_{0}~,\quad \textnormal{pour}~ 1\leq i\leq d-1~,
\end{equation*}
hence we deduce:
\begin{equation*}
\dfrac{1}{k_{0}^{d}}\leq \dfrac{\varphi\left( x\right) }{\varphi\left( y\right) }%
\leq k_{0}^{d}~.
\end{equation*}
Then, noticing that $k_{0}\geq 1$ and taking $k=k_{0}^{D}$~, with $D$ the number
of edges of the subgraph $K$, we obtain
\begin{equation*}
\dfrac{1}{k}\leq \dfrac{\varphi\left( x\right) }{\varphi\left( y\right) }%
\leq k.
\end{equation*}
\end{enumerate}
\end{demorefharnak}

\begin{demorefharmo}
Assume that $\langle Pf,f\rangle \ > 0$ for any function $f\in C_{0
}\left( V\right)\setminus\lbrace0\rbrace$~. Let $x_{0}$ a fixed vertex in $%
V$, which we take as an ``origin''. Consider for $n \geq 1$,  the subgraph $%
G_{n}$ of $G$, such that the set of its vertices is the ball centered in $x_{0}
$ with a radius $n$~ denoted by ${\mathcal{B}}_{n}~,$
\begin{equation*}
{\mathcal{B}}_{n}=\lbrace x\in V;d\left( x_{0},x\right)\leq n\rbrace
\end{equation*}
where $d\left( x,y\right)$ is the combinatorial distance between the vertices $x$
and $y$ in $V$~, which is the number of the edges of the shortest path
connecting $x$ to $y$~. The ball ${\mathcal{B}}_{n}$ is connected and we apply
 Lemma 4.2~, taking it as $K$~, and choosing as function
 $u$ the constant function $1$ on $\partial {\mathcal{B}}_{n}~.$
\newline
We proceed on three steps:
\begin{itemize}
\item {First step:} There exists a function $\psi_{n}\in C_{0}\left(
V\right)$ satisfying $P\psi_{n} \equiv0$~, and such that $\psi_{n}\ >0 $ in the
interior of ${\mathcal{B}}_{n }$ and constant  $1$ on $%
\partial {\mathcal{B}}_{n}$~. Then we consider the function $\Phi_{n }\in
C_{0}\left( V\right)$ given by:
\begin{equation*}
\Phi_{n }\left( x\right)=\dfrac{\psi_{n}\left( x\right) }{\psi_{n}\left(
x_{0}\right)~.}
\end{equation*}
It satisfies the four following conditions:
\begin{enumerate}
\item[\textit{i}.] $\Phi_{n }\left( x_{0}\right) =1$~.

\item[\textit{ii}.] $P\Phi_{n }\equiv0$ in the interior of ${\mathcal{B}}%
_{n}$~.

\item[\textit{iii}.] $\Phi_{n }\upharpoonright \partial {\mathcal{B}}%
_{n}\equiv\dfrac{1} {\psi_{n}\left( x_{0}\right) }$ \; strictly
positive constant.

\item[\textit{iv}.] $\Phi_{n} \ >0$ on ${\mathcal{B}}_{n}~.$
\end{enumerate}

\item {Second step:} Let $x$ a vertex in $V$~, and $n_{0}$ a fixed integer
such that $x$ is in the interior of ${\mathcal{B}}_{n_{0}}$~.

Then for any $n\geq
n_{0}$~, we have: ${\mathcal{B}}_{n_{0}}\subseteq {\mathcal{B}}_{n}$~. Furthermore $%
\Phi_{n}$ is strictly positive ${\mathcal{B}}_{n_{0}}$ and $P-$harmonic
in the interior of ${\mathcal{B}}_{n_{0}}$~. Then from
Lemma 4.1~, there exists a constant $k_{n_{0}}\ > 0 $ such that

\begin{equation*}
\dfrac{1}{k_{n_{0}}}\leq \dfrac{\Phi_{n}\left( x\right) } {\Phi_{n}\left(
x_{0}\right) }\leq k_{n_{0}}~.
\end{equation*}

As $\Phi_{n }\left( x_{0}\right) =1$~, we obtain:

\begin{equation*}
\dfrac{1}{k_{n_{0}}}\leq \Phi_{n}\left( x\right) \leq k_{n_{0}}~.
\end{equation*}

It follows that the set $\lbrace \Phi_{n}\left( x\right) \rbrace_{n\geq
n_{0}}$ is included in the segment $\left[ \dfrac{1}{k_{n_{0}}},k_{n_{0} }%
\right]$~.

\item {Third step:} Let us consider the subset $C$ of ${\mathbb{R}}^{V}$ defined by:
\begin{equation*}
C=\prod_{x\in V} \left[ \dfrac{1}{k_{n_{0}}},k_{n_{0}} \right]~.
\end{equation*}
The sequence $\left( \Phi_{n}\right)_{n\geq n_{0}}$ is in the
compact $C$~, so it has a convergent subsequence $\left( \Phi_{h\left(
n\right) }\right)_{n\geq n_{0}}$ for the topology of ${\mathbb{R}}^{V}$
to a function $\Phi$ satisfying in particular the two following conditions:

\begin{enumerate}

\item[\textit{i}.] $\Phi$ is strictly positive on $V$~, since $%
\Phi\left( x\right)\in \left[ \dfrac{1}{k_{n_{0}}},k_{n_{0}} \right]$~,
for any vertex $x$ in $V$~.

\item[\textit{ii}.] $P\Phi\equiv 0$ on $V$~, since $\underset{%
n\rightarrow\infty}{lim} P\Phi_{h\left( n\right) }\left(x \right)=P\Phi
\left( x\right)$~, for any vertex $x$ in $V$~.

\end{enumerate}
\end{itemize}
\end{demorefharmo}

The function $\Phi$ given by Theorem 4.1 is used to build a unitary transform in
Theorem 5.1~.
%%%%%%%%%%%%%%%%%%%%%%%%%%%%%%%%%%%%%%%%%%%%%%%%%%%%%%%%%%%%%%%%%%%%%%%%%%%%%%%%%%%%%%%%%%%%%%%%%
\section{Any positive Schr\"odinger operator is unitary equivalent to a Laplacian}

We prove that under a positivity condition, a Schr\"odinger operator
is unitary equivalent to a Laplacian
$\Delta_{\omega,c}~.$

\begin{theo}
Let $P$ a Schr\"odinger operator on a graph $G~.$
We assume that
\begin{equation*}
\langle Pf,f\rangle_{l^{2}} > 0
\end{equation*}
for any function $f\in C_0\left( V\right)\setminus\lbrace0\rbrace~.$
\\Then there exist a weight function
 $\omega:V\longrightarrow {\mathbb{R}}_{+}^{\star}$ on $V$ and
a conductance $c:E\longrightarrow {\mathbb{R}}_{+}^{\star}$
on $E$ such that the operator $P$ is unitary equivalent to the
weighted graph Laplacian $\Delta_{\omega,c}$ on $G$~.
\end{theo}

\begin{demo}
We will use a function $\Phi$ which is $P-$harmonique and strictly positive, given by
Theorem 4.1.
Let $P =\Delta_{1,a}+W$ a Schr\"odinger operator
satisfying the hypothesis of the Theorem. By Theorem 4.1~,
 there exists a strictly positive $P-$harmonique function $\Phi$ on $V.$
 Then we obtain:
 $$W=-\dfrac{\Delta\Phi}{\Phi}~.$$
Let us set  $\omega=\Phi$ and for any $g\in l^{2}\left(V \right)~$,
 $f=\dfrac{g}{\Phi}~.$
\newline
We will prove that $\left\langle Pg,g\right\rangle_{l^2}=\left\langle
\Delta_{\omega,c}f,f\right\rangle_{l^2_{\omega}}~.$\\
Let us compute
\begin{align}
\left\langle Pg,g\right\rangle_{l^2}&= \left\langle \Delta \left(
f\Phi\right) +Wf\Phi,f\Phi\right\rangle_{l^2}  \notag \\
&=\left\langle \Delta \left(
f\Phi\right)-f\Delta\Phi,f\Phi\right\rangle_{l^2}  \notag \\
&=\sum_{x\in V}\left[ \sum_{y\sim x} a_{ x,y} \left(f\left(
x\right)\Phi\left( x\right)- f\left( y\right)\Phi\left( y\right)
\right)\right.  \notag \\
&-\left. f\left( x\right) \left( \sum_{y\sim x} a_{ x,y} \left[
\Phi\left( x\right)-\Phi\left(y \right) \right] \right) f\left(
x\right)\Phi\left( x\right) \right]  \notag \\
&=\sum_{x\in V}f\left( x\right)\Phi\left( x\right) \sum_{y\sim x} a_{
x,y} \Phi\left( y\right)\left[ f\left( x\right)-f\left(y \right) %
\right]  \notag \\
&=\sum_{x\in V}\Phi^{2}\left( x\right)f\left( x\right) \dfrac{1}{%
\Phi^2\left( x\right)} \sum_{y\sim x} a_{ x,y}\Phi\left(
x\right)\Phi\left( y\right) \left[ f\left( x\right)-f\left(y \right) \right]~.
\notag
\end{align}
Setting
\begin{equation*}
c_{ x,y}= a_{x,y}\Phi\left( x\right)\Phi\left( y\right)~,
\end{equation*}
 we deduce:
\begin{equation*}
\left\langle Pg,g\right\rangle_{l^2}= \left\langle
\Delta_{\omega,c}f,f\right\rangle_{l^2_{\omega}}.
\end{equation*}
So
\begin{equation*}
P=U^{-1}~\Delta_{\omega,c}~U~,
\end{equation*}
where $U:l^2(V) \longrightarrow l^2_{\omega}(V)$ is given by
\begin{equation*}
U\left( g\right)=\dfrac{g}{\Phi}~.
\end{equation*}
Thus $P$ is unitary equivalent to the Laplacian $\Delta_{\omega,c}$,
with the weight $\omega \equiv \Phi$ and the conductance $c$ given by
\begin{equation*}
c_{x,y}=a_{ x,y }\Phi(x)\Phi(y)~.
\end{equation*}
\end{demo}
%%%%%%%%%%%%%%%%%%%%%%%%%%%%%%%%%%%%%%%%%%%%%%%%%%%%%%%%%%%%%%%%%%%%%%%%%%%%%%%%%%%%%%%%%%%%%%%%%%%%%%%%%%%%%%%%%%%%%%%%%%%%%%%%%%%%%%%%
\section{Metrically complete graphs}

We adapt to graphs the G. and I. Nenciu's method (see [15])
 in the proof of Theorem 6.1~, using Agmon's estimates which are
given by the following technical lemma.

\begin{lemm}
 Let $H =\Delta_{1,a}+W$ a Schr\"odinger operator on
$G,~~\lambda$ a real number and $v\in C\left(V \right)$~.
We assume that $v$ is a solution of the equation:
\begin{equation}
\left(H-\lambda\right)\left( v\right) =0 ~.
\end{equation}
Then for any $f\in C_{0}\left(V\right)$~, we have:
\begin{align}
\left\langle fv,\left( H-\lambda\right) \left( fv\right) \right\rangle_{l^{2}}
&=\sum_{\{x,y \} \in E}a_{ x,y}v\left(x
\right)v\left( y\right) \left[ f\left(x \right)-f\left( y\right) \right]^{2}
\notag \\
&=\dfrac{1}{2}\sum_{x\in V}v\left( x\right) \sum _{y\sim x} a_{
x,y}v\left( y\right) \left[ f\left(x \right)-f\left( y\right) \right]%
^{2}~.  \notag
\end{align}
\end{lemm}

\begin{demo}
Let us assume that: $\left( H-\lambda\right)\left( v\right)=0$~, ie.
for any vertex $x\in V$~,
\begin{equation*}
\sum _{y\sim x} a_{ x,y}\left( v\left( x\right) -v\left(
y\right)\right) + W\left( x\right)v\left( x\right)=\lambda v\left( x\right)
\end{equation*}
Let us compute $S=\left\langle fv,\left( H-\lambda\right) \left( fv\right)
\right\rangle_{l^{2}}$
\begin{align}
S &= \sum_{x\in V}f\left( x\right)v\left( x\right) \left[\left( H-\lambda
\right)\left( fv\right)\right]\left(x \right)  \notag \\
&=\sum_{x\in V}f\left( x\right)v\left( x\right) W\left(x\right)f\left( x\right)v\left( x\right) -
\lambda f\left( x\right)v\left(x\right) \notag \\
&+ \sum_{x\in V} \sum _{y\sim x} a_{ x,y} \left[f\left(
x\right)v\left( x\right)-f\left( y\right)v\left( y\right) \right]~.
\notag
\end{align}
And by the assumption on $v$~, we have:
\begin{equation*}
\lambda f\left( x\right)v\left( x\right)-W\left( x\right)f\left(
x\right)v\left( x\right) = \sum _{y\sim x} a_{ x,y} f\left(
x\right)\left[ v\left( x\right)-v\left( y\right)\right]~.
\end{equation*}
Then, replacing in the precedent expression, we find:

\begin{align}
S &=\sum _{x\in V}f\left( x\right)v\left( x\right)\sum _{y\sim x} a_{
x,y}v\left( y\right) \left[ f\left(x \right)-f\left( y\right) \right]
\notag \\
&=\sum _{x\in V}\sum _{y\sim x} a_{ x,y}v\left(x \right)
v\left( y\right) \left[ f^{2}\left(x \right)-f\left( x\right) f\left(
y\right) \right] \notag
\end{align}
As $a_{x,y }=a_{ y,x }~,$
the expression becomes:
\begin{equation*}
S =\sum_{\{x,y \} \in E} a_{ x,y}v\left(x
\right)v\left( y\right) \left[ f^{2}\left( x\right) -f\left( x\right)
f\left( y\right) + f^{2}\left( y\right) -f\left( x\right) f\left( y\right)%
\right]
\end{equation*}
Finally:
\begin{align}
\left\langle fv,\left( H-\lambda\right) \left( fv\right) \right\rangle_{l^{2}}
&=\sum_{\{x,y \} \in E}a_{ x,y}v\left(x
\right)v\left( y\right) \left[ f\left(x \right)-f\left( y\right) \right]^{2}
\notag \\
&=\dfrac{1}{2}\sum_{x\in V}v\left( x\right) \sum _{y\sim x} a_{
x,y}v\left( y\right) \left[ f\left(x \right)-f\left( y\right) \right]%
^{2}~.  \notag
\end{align}
\end{demo}

\begin{defi}
A graph $G$ is called with bounded degree if there exists an integer $N$ such
that for any $x\in V$ we have: $\sharp \left\lbrace y\in V;~y\sim x\right\rbrace\leq N~.$
\end{defi}

\begin{defi}
Let $a$ a strictly positive function on the edges of the graph $G$~. We define
the $a-$weighted distance on $G$~, which we denote by $\delta_{a}$~:
\begin{equation*}
\delta_{a} \left( x,y\right)=\min _{\gamma\in \Gamma_{x,y}}L\left( \gamma
\right)
\end{equation*}
where $\Gamma_{x,y} $ is the set of all  edge paths $\gamma:x_{1}=x,x_2,$...$x_{n}=y$~,
linking the vertex $x$ to the vertex $y~; $
and
\begin{equation*}
L\left( \gamma\right)=\underset{1\leq i\leq n}{\sum} \dfrac{1}{\sqrt{a_{x_{i}x_{i+1}}}}
\end{equation*}
the length of the edge path $\gamma$~.
\end{defi}

\begin{theo}
 Let $H=\Delta_{1,a}+W$ a Schr\"odinger operator on
an infinite graph $G$ with bounded degree and such that
metric associated to the distance $\delta_{a}$ is complete. We assume
that there exists a real number $k$ so that
\begin{equation*}
\left\langle Hg,g \right\rangle_{l^{2}} \geq k\Vert g\Vert^{2}_{l^{2}}
\end{equation*}
for any $g\in C_{0 }\left( V\right)$.
Then the operator $H$~, with domain $C_{0 }\left( V\right)$, is essentially self-adjoint.
\end{theo}

\begin{demo}
Let $\lambda \ < k-1$~, we would prove that if $v\in l^{2}\left( V\right)$ and
satisfies the equation
\begin{equation*}
Hv=\lambda v~,
\end{equation*}
 then $v$  vanishes.%
\newline
We set $R\ >0$ and a vertex $x_{0}$ as the origin. Let as denote:
\begin{equation*}
B_{R}=\left\lbrace x\in V;\delta_{a} \left( x_{0},x\right )\leq
R\right\rbrace
\end{equation*}
 the ball centered in $x_{0}$ and with radius $R$ for the distance
$\delta_{a}$.\newline
Let us consider the function $f$ defined on $V$ by
\begin{equation*}
f\left( x\right) =\min\left( 1,\delta_{a}\left(x, V\setminus
B_{R+1}\right)\right)~.
\end{equation*}
Hence we have:
\begin{equation*}
f\upharpoonright B_{R}\equiv 1,\quad f\upharpoonright V\setminus
B_{R+1}\equiv0~,  \quad f\left( B_{R+1}\setminus B_{R} \right)\subseteq \left[
0,1\right]
\end{equation*}
The support of $f$ is in the bounded ball $B_{R+1}$ which is finite dy the completeness of the
metric associated to  the distance $\delta_{a}$.\newline
By the assumption on $H$ and the finiteness of the support of $fv$ in $%
B_{R+1}$~, we find the following inequalities:
\begin{equation*}
\left\langle fv,\left( H-\lambda\right) \left( fv\right) \right\rangle_{l^{2}} \geq
\left( k-\lambda \right) \sum_{x\in B_{R+1}}\left( fv\right) ^2\left( x\right)
\geq\sum _{x\in B_{R}} v^2\left( x\right)
\end{equation*}
On the other hand, using Lemma 6.1~, we find:
\begin{align}
\left\langle fv,\left( H-\lambda \right) \left( fv\right) \right\rangle_{l^{2}}&=
\dfrac{1}{2}\sum _{x\in V}\sum_{y\sim x}a_{ x,y}
v\left( x\right)v\left(y \right) \left[f\left( x\right)-f\left( y\right) %
\right]^2  \notag \\
&\leq \dfrac{1}{2}\sum _{x\in V}\sum_{y\sim x}a_{
x,y} v^{2}\left( x\right) \left[f\left( x\right)-f\left(
y\right) \right]^2 \notag
\end{align}
using the fact that $a_{x,y}=a_{y,x}$ and that:
\begin{equation*}
v\left( x\right)v\left(y \right) \leq\dfrac{1}{2}\left(v^{2}\left(
x\right)+v^{2}\left( y\right) \right)~.
\end{equation*}
As every restriction of $f$ to $B_{R}$ and to $V\setminus B_{R+1}$
are constant functions, then the preceding inequality becomes:
\begin{align}
\left\langle fv,\left( H-\lambda\right) \left( fv\right) \right\rangle_{l^{2}}&\leq
\dfrac{1}{2}\sum _{x\in B_{R+1}\setminus B_{R}}\sum_{y\sim x}a_{
x,y} v^2\left( x\right) \left[f\left( x\right)-f\left( y\right) %
\right]^2  \notag \\
&\leq\dfrac{1}{2}\sum _{x\in B_{R+1}\setminus B_{R}}\sum_{y\sim
x}a_{ x,y} v^{2}\left( x\right)\left(
\delta_{a}\left( x,y\right) \right)^{2}  \notag
\end{align}
The last inequality is obtained by the fact that $f$ is a $1-$Lipschitz
 continuous function since it is the minimum of two
 $1-$Lipschitz continuous functions.\newline
Since the distance $\delta_{a}$ satisfies the condition:
\begin{equation*}
\delta_{a}\left( x,y\right)\leq\frac{1}{\sqrt{a_{
x,y}}}
\end{equation*}
if $\left\lbrace x,y\right\rbrace$ is an edge, and the degree of $G$
is bounded by $N$~, we have:
\begin{equation*}
\left\langle fv,\left( H-\lambda\right) \left( fv\right) \right\rangle_{l^{2}} \leq%
\dfrac{1}{2}N\sum_{x\in B_{R+1}\setminus B_{R}}v^2\left( x\right)
\end{equation*}
Hence, for any $R\ >0$~, we find:
\begin{equation*}
\sum _{x\in B_{R}}v^2\left( x\right)\leq\left\langle fv,\left(
H-\lambda\right) \left( fv\right) \right\rangle_{l^{2}} \leq\dfrac{1}{2}N\sum_{x\in
B_{R+1}\setminus B_{R}}v^2\left( x\right)
\end{equation*}
And since $v\in l^2\left( V\right)$~, making $R$ tend to $\infty$~,
 it results that:
\begin{equation*}
\lim_{R\rightarrow\infty}\sum _{x\in B_{R+1}\setminus B_{R}}v^2\left(
x\right) =0.
\end{equation*}
It follows $\Vert v\Vert^2_{l^2} =0$~.
\end{demo}

\begin{theo}
Let $G$ an infinite graph with bounded degree and which is weighted
by $\omega$ on $V$ and a conductance $c$ on $E$.
We assume that the metric associated to the distance
$\delta_{a}$ is complete, where $a$ is the function given by
\begin{equation*}
a_{ x,y}= \dfrac{%
c_{ x,y}}{\omega_{x}\omega_{y}}~.
\end{equation*}
Then the Laplacian $\Delta _{\omega,c~,}$ with domain $C_0(V)$, is
essentially self-adjoint.
\end{theo}

\begin{demo}
By Lemma 2.1, the operator $\Delta _{\omega,c}$ is
unitary equivalent to the Schr\"odinger operator  $%
H=\Delta_{1,a}+W$~, where
\begin{equation*}
a_{x,y }= \dfrac{c_{ x,y }%
}{\omega_{x}\omega_{y}}
\end{equation*}
and
\begin{equation*}
W\left( x\right)= \dfrac{1}{\omega_{x}^{2}}~\sum_{y\sim x}
c_{x,y} \left( 1-\dfrac{\omega_{x}}{\omega_{y}}\right).
\end{equation*}
And we apply  Theorem 6.1 for the operator $H$ which
satisfies the assumptions, since $\Delta_{\omega,c}$ is positive and
\begin{equation*}
\left\langle Hg,g \right\rangle_{l^2}= \left\langle \Delta_{\omega,c}%
\dfrac{g}{\omega} ,\dfrac{g}{\omega} \right\rangle_{l_{\omega}^2}
\end{equation*}
in the proof of Theorem 5.1~.
\end{demo}

\begin{rem}
 Theorem 6.1 is not particular case  of Theorem 3.2.
 In fact in Theorem 6.1, the potential $W$ is
not necessarily  bounded from below.
For example let $G$ the graph  such that
$V= {\mathbb{N}} \setminus\ \left\lbrace 0,1\right\rbrace $
 and $n\sim n+1$ for any $n$~. We assume that $G$ is
weighted by the vertex weight $\omega_n=\dfrac{1}{n\ \log n}$ and by
the constant edge conductance $c_n=1$~. The distance $\delta_{a}$ is given by:
\begin{equation*}
\delta_{a}( n_0,n )=\sum_{n_0\leq k\leq n} \dfrac{1}{\sqrt{a_{k,k+1}}}
\end{equation*}
but
\begin{equation*}
\dfrac{1}{\sqrt{a_{k,k+1}}} = \dfrac{1}{\sqrt{k(k+1)\log k\log(k+1)}}%
\underset{\infty}{\sim}\dfrac{1}{k\log k}~,
\end{equation*}
then $$\delta_{a}\ ( n_0,n\ )\underset{n\rightarrow\infty}{%
\longrightarrow}\infty~,$$
and the associated metric is complete.\\
Furthermore, setting $$H=\Delta_{1,a}+W~,$$
we have: $$\left\langle Hg,g
\right\rangle_{l^2}= \left\langle \Delta_{\omega,c}\dfrac{g}{\omega},%
\dfrac{g}{\omega} \right\rangle_{l_{\omega}^2}\geq0~,$$
for any $g\in C_0(V)~.$\\

While the potential $W$ is not bounded from below, since we obtain, after calculation:
\begin{equation*}
W\left( n\right)=2n^2\log^2n-n\log n\left[(n+1)\log(n+1)+ (n-1)\log(n-1)%
\right]\underset{\infty}{\sim}-\log n
\end{equation*}
which goes to $-\infty~.$
\end{rem}

\begin{rem}
In the Example of Remark 6.1~, the choice of the weight according to $log$ is crucial.
 In fact, setting power functions, we can not have at the same time
the metric $\delta_{a}~$ complete and
the potential $W$ not bounded from below.\\
For example let  $G$ the graph such that $V= {\mathbb{N}} \setminus\ \left\lbrace
0,1\right\rbrace $ and $n\sim n+1$ for any $n~,$ and assume $G$ weighted by the vertex weight
$\omega_n=\dfrac{1}{n^{\alpha}}$ and the edge conductance
$c_n=\dfrac{1}{n^{\beta}}.$
The distance $\delta_{a}~$ is given by:
\begin{equation*}
   \delta_{a}( n_0,n )=\sum_{n_0\leq k\leq n}\left(\dfrac{k^{\beta}}{(k^{\alpha}(k+1)^{\alpha})}\right)^{\frac{1}{2}}~.
\end{equation*}
And an easy calculation show that
\begin{equation*}
   \delta_{a}~\textnormal{is complete  if and only if}~\alpha-\dfrac{1}{2}~\beta \leq 1~.
\end{equation*}
We obtain easily for the potential:
\begin{equation*}
   W_{n}\sim -\alpha(\alpha -\beta-1)n^{2\alpha-\beta-2}~.
\end{equation*}
Hence for $\alpha-\dfrac{1}{2}~\beta \leq 1~,$ the potential $W$
is bounded from below.
\end{rem}

\begin{rem}
The completeness of the metric $\delta_{a}$ is not a necessary condition
to essential self-adjointeness of the Laplacian $\Delta_{\omega,c}$ (or
the Schr\"{o}dinger operator $\Delta_{1,a}+W~$.)\\
 In fact, let $G$ the graph $\N$. For any edge conductance
$a$, the Laplacian $\Delta_{1,a}~$
is essentially self-adjoint by Theorem 3.1~.
While the metric given by the distance $\delta_{a}~$
is not necessarily  complete: for example when
$a_{n}=(n+1)^{-2-\varepsilon}~,$ for an $\varepsilon > 0~.$
\newline

We give in the paper [5] some increasing conditions of the
potential insuring the essential self-adjointness of a Schr\"{o}dinger operator
 on metrically non complete graphs.
\end{rem}

\begin{rem}\label{K-L}
  Theorem 6.2 and Keller-Lenz's Theorem are not deduced one from another. For example, let $G$
 the graph $\N$ . We consider the two following cases:
 \begin{enumerate}
   \item  Let us choose the vertex weight $\go_{n}=\frac{1}{n+1}$ and the constant edge conductance
 $c_{n}=1$. For a fixed vertex $n_{0}$, we have
 \begin{equation*}
   \gd_{a }(n_{0},n)=\sum_{n_{0}\leq k\leq n}\frac{1}{\sqrt{k(k+1)}}\underset{n\rightarrow\infty}{%
   \longrightarrow}\infty~,
\end{equation*}
  So the associated metric is complete. Hence by  Theorem 6.2, we conclude that the
  Laplacian $\Delta_{\go,1}$ is essentially self-adjoint.\\
   But the sum  $\sum_{n} \go _{n}^{2}$~~ is finite, so the assumption $(A)$ is not satisfied.
    Hence we can not conclude by  Keller-Lenz's Theorem.
\item Let us choose now the vertex weight $\go_{n}=\frac{1}{\sqrt{n+1}}$ and the edge conductance
 $c_{n}=n^{2}$. Then
 \begin{equation*}
   \gd_{a }(n_{0},n)=\sum_{n_{0}\leq k\leq n} \frac{1}{k((k+1)(k+2))^{\frac{1}{4}}},
\end{equation*}
which is convergent. Then the associated metric is non complete, and we can not apply Theorem 6.2 here. While
Keller-Lenz's Theorem can be applied as $\sum_{n} \go _{n}^{2}$~~ is not finite.
 \end{enumerate}
    \end{rem}
%%%%%%%%%%%%%%%%%%%%%%%%%%%%%%%%%%%%%%%%%%%%%%%%%%%%%%%%%%%%%%%%%%%%%%%%%%%%%%%%%%%%%%%%%%%%%%%%%%%%%%%%%%%%%%%%%%%%%%%%%%%%%%%%%%%%%%%%%%

\end{document}